\documentclass[11pt,a4paper]{article}
\usepackage[margin=25mm]{geometry}
\usepackage{parskip}
\usepackage[T1]{fontenc}
\usepackage{mathpazo}
\RequirePackage[english]{babel}
\usepackage{amsmath,amssymb,amsthm}
\usepackage{tikz}
\usetikzlibrary{arrows.meta}
\usepackage{placeins}
\usepackage[shortlabels]{enumitem}
\setlist{itemsep=0ex,topsep=0ex,parsep=0.4ex}
\newtheorem{theorem}{Theorem}
\newtheorem{lemma}[theorem]{Lemma}
\newtheorem{corollary}[theorem]{Corollary}
\newtheorem{conjecture}[theorem]{Conjecture}
\newtheorem{prop}[theorem]{Proposition}
\newtheorem{Question}[theorem]{Question}

\newtheorem{claim}{Claim}
\newenvironment{claimproof} {\begin{proof}}{\end{proof}}

\makeatletter
\renewenvironment{proof}[1][\proofname]
{\par\pushQED{\qed}
  \normalfont\topsep6\p@\@plus6\p@\relax\trivlist
  \item[\hskip\labelsep\bfseries#1\@addpunct{.}]
  \ignorespaces}
{\popQED\endtrivlist\@endpefalse}
\makeatother

\DeclareMathOperator{\diam}{diam}
\DeclareMathOperator{\dist}{dist}

\DeclareMathOperator{\tparsymbol}{\mathit{t}}
\newcommand{\tpar}{\mathord{\tparsymbol}}
\newcommand{\C}{\mathcal R}
\newcommand{\floor}[1]{\left\lfloor #1\right\rfloor}

\definecolor{colourAone}{HTML}{F0C3C3}
\definecolor{colourAtwo}{HTML}{F3DDA0}
\definecolor{colourAthree}{HTML}{B9DDD6}
\definecolor{colourBone}{HTML}{BDC3FF}
\definecolor{colourBtwo}{HTML}{D9C5E5}
\definecolor{colourBthree}{HTML}{C7DCAF}
\tikzset{
  coloured vertex/.style={circle,draw=black!75,line width=.6pt,
    minimum size=9.5mm,inner sep=1pt,align=center,font=\small},
  graph edge/.style={draw=black!65,line width=.65pt},
  recolouring arrow/.style={-{Stealth[length=2.5mm]},line width=.85pt}
}
\newcommand{\colourvertex}[5]{%
  \node[coloured vertex,fill=#5] (#1) at #2
    {\shortstack{\scriptsize $#3$\\[-1pt]$#4$}};
}

\newcommand{\Yvertexsize}{7.7mm} 
\tikzset{
  square vertices/.style={
    coloured vertex/.append style={
      rectangle,
      minimum size=\Yvertexsize,
      inner sep=1pt
    }
  }
}

\usepackage{xcolor}
\definecolor{cblue}{RGB}{0,70,140}
\definecolor{cgreen}{RGB}{100,140,0}
\definecolor{cred}{RGB}{190,10,50}

\usepackage[hyphens]{url}
\usepackage[colorlinks,pagebackref]{hyperref}
\hypersetup{citecolor=cgreen,linkcolor=cblue,urlcolor=cblue}
\usepackage[capitalise,nameinlink,noabbrev,compress]{cleveref}

\renewcommand*{\backref}[1]{}
\renewcommand*{\backrefalt}[4]{%
    \ifcase #1 Not cited.%
    \or $\uparrow$#2%
    \else $\uparrow$#2%
    \fi%
}
\let\oldbibliography\bibliography
\renewcommand{\bibliography}[1]{%
  {%
    \hypersetup{linkcolor=cred}%
    \bibliographystyle{bibstyle}%
    \oldbibliography{#1}%
  }%
}

\newcommand{\affiliation}{\footnote}
\makeatletter
\def\@fnsymbol#1{\ensuremath{\ifcase#1\or *\or \dagger\or \ddagger\or \mathsection\or \|\or **\or \dagger\dagger \or \ddagger\ddagger \else\@ctrerr\fi}}
\makeatother

\author{
W. Cames van Batenburg\affiliation{D\'epartement d'Informatique, Universit\'e libre de Bruxelles, Belgium (\textsf{\href{mailto:w.p.s.camesvanbatenburg@gmail.com}{w.p.s.camesvanbatenburg@gmail.com}}). Supported by the Belgian National Fund for Scientific Research (FNRS).}}

\title{Large diameter for $(\Delta+2)$-colouring reconfiguration}
\begin{document}
\maketitle

\begin{abstract}
Given a graph $G$ and an integer $k$, the reconfiguration graph $\C_{k}(G)$ has the proper $k$-colourings of $G$ for vertices, two of which are adjacent if and only if they differ on precisely one vertex of $G$. In $2022$,
Cambie, Cranston and the author proved that the diameter of $\C_{\Delta+2}(G)$ is at most $n+2\mu$, and put forth a collection of conjectures that in its strongest and most beautiful form would imply that the diameter of $\C_{\Delta+2}(G)$ is precisely equal to $n + \mu$, for every graph $G$ with $n$ vertices, matching number $\mu$ and maximum degree $\Delta$. Here we disprove most of those conjectures. Our constructions have reconfiguration diameter at least $n+\frac{3}{2} \mu -2$. We also provide counterexamples that are $\Delta$-regular, we show that every graph has reconfiguration diameter at least $n+\frac{1}{2} \mu$, and we advertise an extension to hypergraphs. 
It is now much more mysterious 
how the diameter behaves,
so we propose a number of open problems.
\end{abstract}
\section{Introduction}
Given a set of objects $\mathcal S$ and a prescribed collection of reversible moves, one can define an associated \emph{reconfiguration graph} $\mathcal{R}$, with vertex set $\mathcal S$, and two distinct objects adjacent if and only if one can be obtained from the other by a single allowed move. Typically, the moves are chosen to reflect a small local modification of the object, so that the distance in $\mathcal{R}$ between $s,t\in \mathcal{S}$ is a good measure for how closely related $s$ and $t$ are.  More precisely, reconfiguration graphs form a model for how easily two objects can be converted into each other via a sequence of (small) changes, while staying inside the space of relevant objects.

In the classic example, $\mathcal{S}$ is the set of configurations of a Rubik's cube, and the allowed moves are the turns of the cube. Then the diameter of $\C$ is the maximum number of turns required to convert any two given configurations into each other. In particular, this upper bounds the maximum number of moves needed to solve any Rubik's cube. In more serious situations, $\mathcal{S}$ may represent the possible states of a physical, chemical or biological system. We refer to the surveys~\cite{Nishimura18,Heuvel13} for various other examples.

In one of the most well-studied settings, $\mathcal{S}$ is the set of proper $k$-vertex-colourings of some fixed graph $G$, for some positive integer $k$.
Under the canonical reconfiguration rule, two $k$-colourings $c_1,c_2$ are adjacent if and only if they differ on precisely one vertex of $G$. From now on we only consider this type of reconfiguration graph, and we denote it by $\C_k(G)$.

Connectivity and the diameter of $\C_k(G)$ have been studied
in relation to various graph parameters and structural
restrictions, including
degeneracy~\cite{DyerFlaxmanFriezeVigoda06, BousquetPerarnau16,BousquetHeinrich22}, maximum
average degree~\cite{Feghali21,Cranston22, ZhengLu26}, bipartite graphs \cite{CambieCamesCranston26, AsgarliKrehbielMacLeanZaimi25}, treewidth~\cite{BonamyBousquet13,BonamyBousquet18, BartierBousquetHeinrich21,WangLu26},
excluded minors~\cite{GuLeiLiMaWang26}, planar
graphs~\cite{DvorakFeghali21,BartierBousquetFeghaliHeinrichMoorePierron23,CranstonMahmoud24,Cranston25, DaiLiu26}, and excluded
induced subgraphs (especially for $k$ equal to one plus the chromatic number)~\cite{BonamyJohnsonLignosPatelPaulusma14,BonamyBousquet18, BelavadiCameron25,BelavadiCameronHildred25,LeiMaMiaoShiWang24,BelavadiCameronMerkel24,BelavadiCameronSintiari26}.
In this current work, we focus on the maximum
degree~\cite{jerrum1995very, FeghaliJohnsonPaulusma16,BonamyBousquetPerarnau21, CCC, BousquetFeuilloleyHeinrichRabie25,CambieCamesCranstonHeuvelKang25,DeMeyer26}.

Let $\Delta=\Delta(G)$ and $n=n(G)$ denote the maximum degree and number
of vertices of $G$, respectively. Complete graphs show that $\C_k(G)$
can be empty or disconnected when $k\leq\Delta+1$. However, it is
connected (and hence has finite diameter) whenever $k\geq\Delta+2$. In this range,
Jerrum~\cite{jerrum1995very} in fact gave an argument that implies an exponential upper bound
$2^{O(n)}$ on the diameter. Cereceda~\cite{cereceda2007diss}
subsequently improved this to $O(n\Delta)$, and hence $O(n^2)$.
Finally, Cambie, the author, and Cranston~\cite{CCC}
showed that the diameter is linear in $n$, namely at most $2n-1$.

This is close to optimal: by relabelling colours one obtains two colourings that differ at every vertex, so for every graph $G$ the diameter of $\C_{\Delta+2}(G)$ must be at least $n$. 
Inspired by the fact that $\diam \C_{\Delta+2}(K_n)=\lfloor \tfrac{3n}{2} \rfloor$ for complete graphs $K_n$, and via the intermediate Conjecture~\ref{conj:oldmotivatingconjecture} mentioned at the end of this manuscript, three refined conjectures on the exact value of the diameter were formulated.
Recall that the \emph{matching number} $\mu=\mu(G)$ of a graph $G$ is the maximum number of vertex-disjoint edges and the \emph{vertex cover number} $\tau=\tau(G)$ is the smallest size of a vertex set that contains at least one endpoint of every edge. The inequality  $\mu(G)\le \tau(G)$ holds for every graph, with equality for bipartite graphs due to K\"onig's Theorem. \begin{conjecture}\label{conjectures}\cite{CCC}
For every graph $G$ with $n$ vertices, maximum degree $\Delta$, matching number $\mu$ and vertex cover number $\tau$,
\begin{enumerate}[(i)]
    \item\label{conj1} $\diam\C_{\Delta+2}(G) \le n+ \tau$;
    \item \label{conj2} 
    $\diam\C_{\Delta+2}(G) \le n+ \mu$;
    \item\label{conj3} 
    $\diam\C_{\Delta+2}(G)= n+ \mu$ \text{if} G \text{is regular}.
\end{enumerate}
\end{conjecture}

On the lower bound side, it was also asked in~\cite{CCC} whether $\diam \C_{\Delta+2}(G)\ge n+\mu$ for every graph $G$. This suggested the tantalizingly beautiful possibility of equality $\diam \C_{k}(G)= n+\mu$ for every graph $G$ and $k\ge \Delta+2$.
The following proposition 
gives
a nonexhaustive selection of evidence from the literature that was gathered in support of this, and Conjecture~\ref{conjectures} in particular.

\begin{prop}\label{thm:literature}\cite{CCC,DeMeyer26,Zeven23}
For every graph $G$ with $n$ vertices, maximum degree $\Delta \ge 1$ and matching number $\mu$,
\begin{enumerate}[(i)]
    \item\label{prop1} $\diam\C_{k}(G) = n+ \mu$  \text{ for all} $k\ge 2\Delta+1$;
    \item\label{prop2} $\diam\C_{k}(G) \le n+ 2\mu$ \text{ for all} $k\ge \Delta+2$;
    \item\label{prop3} $\diam\C_{k}(G)\ge n+ \mu$ \text{ for all} $k\ge \Delta+2$  \text{, provided}  $G$ \text{is bipartite};    
    \item\label{prop4} $\diam\C_{k}(G)= n+ \mu$ \text{ for all} $k\ge \Delta+2$  \text{, provided} $\Delta \leq 3$ \text{or} $G$ \text{is complete multipartite or a tree}.
\end{enumerate}  
\end{prop}


\subsection*{Counterexamples}
Conjectures~\ref{conjectures}.\ref{conj1} and~\ref{conjectures}.\ref{conj2} 
are special cases of the \emph{Correspondence Colouring Reconfiguration Conjecture} and the \emph{List Colouring Reconfiguration Conjecture}, respectively, and Conjecture~\ref{conjectures}.\ref{conj3} was called the \emph{Regular Cereceda Conjecture} in~\cite{CCC}. 
We disprove all of them via Theorems~\ref{thm:lower} and~\ref{thm:regular} below.

\begin{theorem}\label{thm:lower}
For every integer $t \ge 2$, there is a graph $G$ with $n(G)=t^2$ vertices, maximum degree $\Delta(G)=2t-2$ and matching number and vertex cover number $\tau(G)=\mu(G)=t$, such that 
$$\diam\C_{\Delta(G)+2}(G) \ge n(G)+ \tfrac{3}{2}\,\mu(G) -2.$$
\end{theorem}


\begin{theorem}\label{thm:regular}
There exist $8$-regular graphs $G$ such that $\diam\C_{10}(G) \ge n(G)+ \tfrac{166}{165}\,\mu(G).$
Moreover, for all even $\Delta \ge 8$ there exists a $\Delta$-regular graph $G$ such that $$\diam\C_{\Delta+2}(G) > n(G)+ \mu(G).$$ 
\end{theorem}

\subsection*{Additional results}
Each graph constructed in the proof of Theorem~\ref{thm:lower} is the $1$-subdivision of a multigraph. Therefore they are bipartite, such that all vertices in one of the two parts have degree $2$. The following converse theorem shows that the factor $3/2$ is the 
asymptotically optimal coefficient of $\mu$ within that class of bipartite graphs.
\begin{theorem}\label{thm:upper}
Let $G$ be a  bipartite graph such that all vertices in one of the parts have degree at most $2$.
If $k= \Delta(G)+2$, then $$\diam\C_{k}(G)\le n(G)+\tfrac{3}{2}\,\mu(G).$$
If $k\ge \Delta(G)+3$, then $$\diam\C_{k}(G)= n(G)+\mu(G).$$
\end{theorem}

We complement this with a simple bound for bipartite graphs in terms of only the number of vertices, which appears to have eluded explicit statement in the literature. It is sharp due to balanced complete bipartite graphs.
\begin{lemma}\label{lem:bipartite}
Let $G$ be a bipartite graph. 
If $k\ge\Delta(G)+2$, then $\diam\C_k(G)\le \tfrac{3}{2}\,n(G).$
\end{lemma}


Finally, we provide a lower bound that is valid for all graphs.
\begin{theorem}\label{thm:generallowerbound}
Let $G$ be a graph. If $k\ge \Delta(G)+2$, then
$$\diam \C_k(G) \ge n(G)+ \tfrac{1}{2} \mu(G).$$
\end{theorem}

In a different direction, we present an extension to hypergraphs $\mathcal{H}$, showing that in that setting  the diameter of the reconfiguration graph is linear in the number of vertices as well. The exact value of the diameter may be more elusive than for graphs, even when the number of colours is large. We refer to Section~\ref{sec:hypergraph} for the relevant definitions and two illustrative corollaries concerning acyclic colourings and partitions into induced forests. 
\begin{lemma}\label{lem:hypergraph_intro}
Let $\mathcal{H}$ be a loopless hypergraph on $n$ vertices and let $k\geq \Delta^*(\mathcal{H})+2$. Then
$$n \leq \diam \C_k(\mathcal{H}) \leq 2n.$$
\end{lemma}

\subsection*{Organization}
The main result  Theorem~\ref{thm:lower} is proved in Section~\ref{sec:lower}. After that, in Section~\ref{sec:regular} we modify the construction to obtain \emph{regular} graphs that certify Theorem~\ref{thm:regular}.  In the remaining sections we prove Lemmas~\ref{lem:bipartite} and~\ref{lem:hypergraph_intro} and stronger forms of Theorem~\ref{thm:upper} and Theorem~\ref{thm:generallowerbound}. 
We end with a number of open problems.

\section{Highly nonregular counterexamples}
\label{sec:lower}
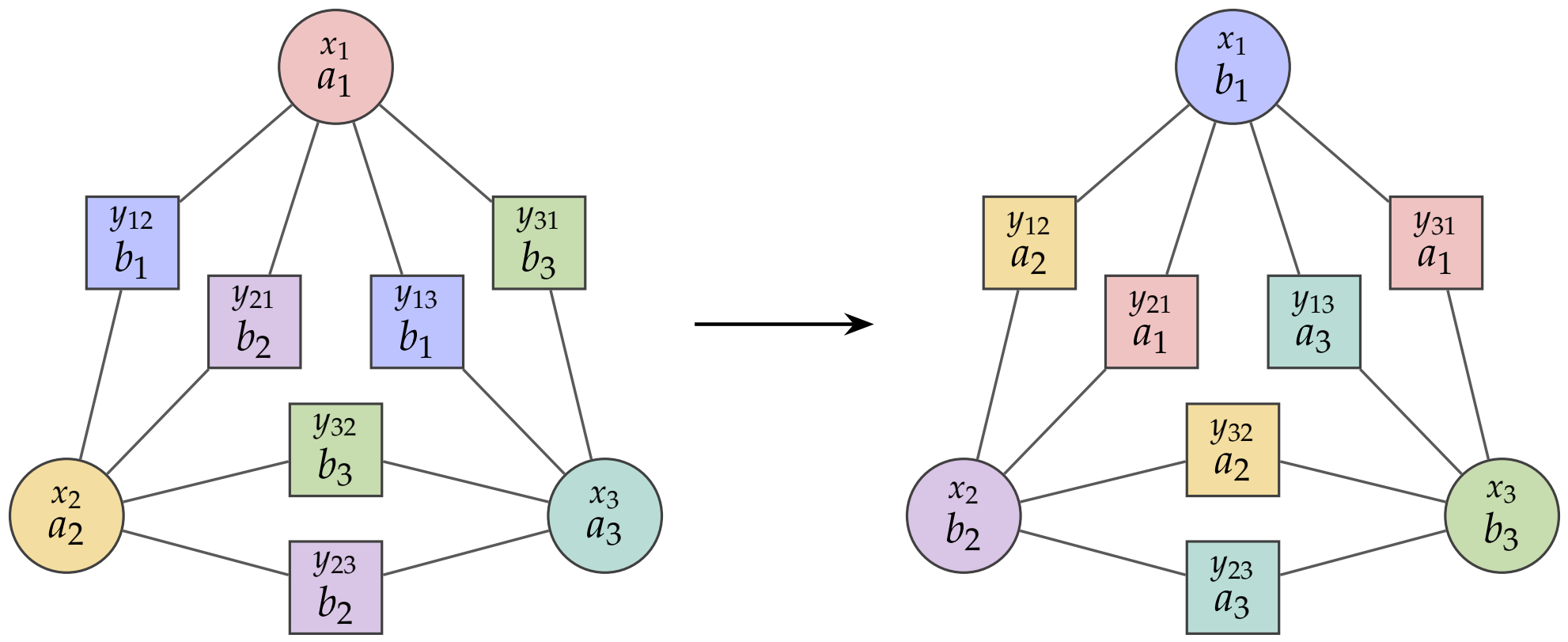
\begin{figure}[htbp]
\centering
\begin{tikzpicture}

\begin{scope}[xshift=0cm]

  \colourvertex{gtAx1}{(0,2.25)}{x_1}{a_1}{colourAone}
  \colourvertex{gtAx2}{(-2.25,-1.5)}{x_2}{a_2}{colourAtwo}
  \colourvertex{gtAx3}{(2.25,-1.5)}{x_3}{a_3}{colourAthree}

  \begin{scope}[square vertices]
    \colourvertex{gtAy12}{(-1.7,0.78)}{y_{12}}{b_1}{colourBone}
    \colourvertex{gtAy21}{(-0.68,0.12)}{y_{21}}{b_2}{colourBtwo}

    \colourvertex{gtAy13}{(0.68,0.12)}{y_{13}}{b_1}{colourBone}
    \colourvertex{gtAy31}{(1.7,0.78)}{y_{31}}{b_3}{colourBthree}

    \colourvertex{gtAy23}{(0,-2.1)}{y_{23}}{b_2}{colourBtwo}
    \colourvertex{gtAy32}{(0,-0.95)}{y_{32}}{b_3}{colourBthree}
  \end{scope}

  \draw[graph edge] (gtAx1)--(gtAy12)--(gtAx2);
  \draw[graph edge] (gtAx1)--(gtAy13)--(gtAx3);
  \draw[graph edge] (gtAx2)--(gtAy21)--(gtAx1);
  \draw[graph edge] (gtAx2)--(gtAy23)--(gtAx3);
  \draw[graph edge] (gtAx3)--(gtAy31)--(gtAx1);
  \draw[graph edge] (gtAx3)--(gtAy32)--(gtAx2);
\end{scope}

\begin{scope}[xshift=7.5cm]
  \colourvertex{gtBx1}{(0,2.25)}{x_1}{b_1}{colourBone}
  \colourvertex{gtBx2}{(-2.25,-1.5)}{x_2}{b_2}{colourBtwo}
  \colourvertex{gtBx3}{(2.25,-1.5)}{x_3}{b_3}{colourBthree}

  \begin{scope}[square vertices]
    \colourvertex{gtBy12}{(-1.7,0.78)}{y_{12}}{a_2}{colourAtwo}
    \colourvertex{gtBy21}{(-0.68,0.12)}{y_{21}}{a_1}{colourAone}

    \colourvertex{gtBy13}{(0.68,0.12)}{y_{13}}{a_3}{colourAthree}
    \colourvertex{gtBy31}{(1.7,0.78)}{y_{31}}{a_1}{colourAone}

    \colourvertex{gtBy23}{(0,-2.1)}{y_{23}}{a_3}{colourAthree}
    \colourvertex{gtBy32}{(0,-0.95)}{y_{32}}{a_2}{colourAtwo}
  \end{scope}

  \draw[graph edge] (gtBx1)--(gtBy12)--(gtBx2);
  \draw[graph edge] (gtBx1)--(gtBy13)--(gtBx3);
  \draw[graph edge] (gtBx2)--(gtBy21)--(gtBx1);
  \draw[graph edge] (gtBx2)--(gtBy23)--(gtBx3);
  \draw[graph edge] (gtBx3)--(gtBy31)--(gtBx1);
  \draw[graph edge] (gtBx3)--(gtBy32)--(gtBx2);
\end{scope}

\draw[recolouring arrow] (3.0,.10)--(4.5,.10);

\end{tikzpicture}
\caption{The graph $G_3$ with the colourings $\alpha$ (left) and $\beta$ (right).}
\label{fig:G3}
\end{figure}

In this section, we prove Theorem~\ref{thm:lower}: for every $\mu\ge 2$ there exists a graph $G_{\mu}$ on $n$ vertices with matching number $\mu$ and maximum degree $\Delta$ such that the reconfiguration graph of proper $(\Delta+2)$-colourings of $G_{\mu}$ has diameter at least $n+ \frac{3 \mu}{2} -2$. The construction is a subdivision of a complete multigraph. See Figure~\ref{fig:G3} for an illustration of $G_3$.

We now give the details.
\begin{proof}[Proof of Theorem~\ref{thm:lower}]
For $\tpar\ge2$, let $G_{\tpar}$ be the bipartite graph with parts
$$X=\{x_1,\ldots,x_{\tpar}\} \ \ \text{ and }  \ \
 Y=\{y_{ij}:1\le i,j\le\tpar,\ i\ne j\},$$
where the neighbourhood of $y_{ij}$  is $\{x_i,x_j\}$, for every $i\neq j$. Thus $G_{\tpar}$ is obtained by subdividing
once every edge of the multigraph formed by doubling each edge of the complete graph
$K_{\tpar}$.
The number of vertices, maximum degree, matching number and vertex cover number of $G_t$ are
$$n(G_{\tpar})=\tpar^2,\qquad \Delta(G_{\tpar})=2{\tpar}-2,\qquad \mu(G_{\tpar})= \tau(G_{\tpar}) =\tpar.$$

Indeed, $X$ is a vertex cover so that $\mu(G_t)\le \tau(G_t)\le t$, and the edges $x_iy_{i,i+1}$ (with indices
modulo $\tpar$) form a matching 
so that $t\le \mu(G_t)$.
We now define two $2t$-colourings $\alpha$ and $\beta$, and we will show that they are at distance at least $n(G_t)+\frac{3}{2} \mu(G_t) -2$. Use $2{\tpar}$ distinct colours $a_1,\ldots,a_{\tpar},b_1,\ldots,b_{\tpar}$,
and define $\alpha$ and $\beta$ by
\begin{equation*}\label{table:colouring_scheme}
\begin{array}{c|cc}
 &\alpha&\beta\\\hline
 x_i&a_i&b_i\\
 y_{ij}&b_i&a_j
\end{array}
\end{equation*}

The restrictions of these colourings to the subgraph induced by $\{x_i,y_{ij},y_{ji},x_j\}$ are shown in Figure~\ref{fig:four-vertices}, and the colourings of $G_3$ are illustrated in Figure~\ref{fig:G3}. 

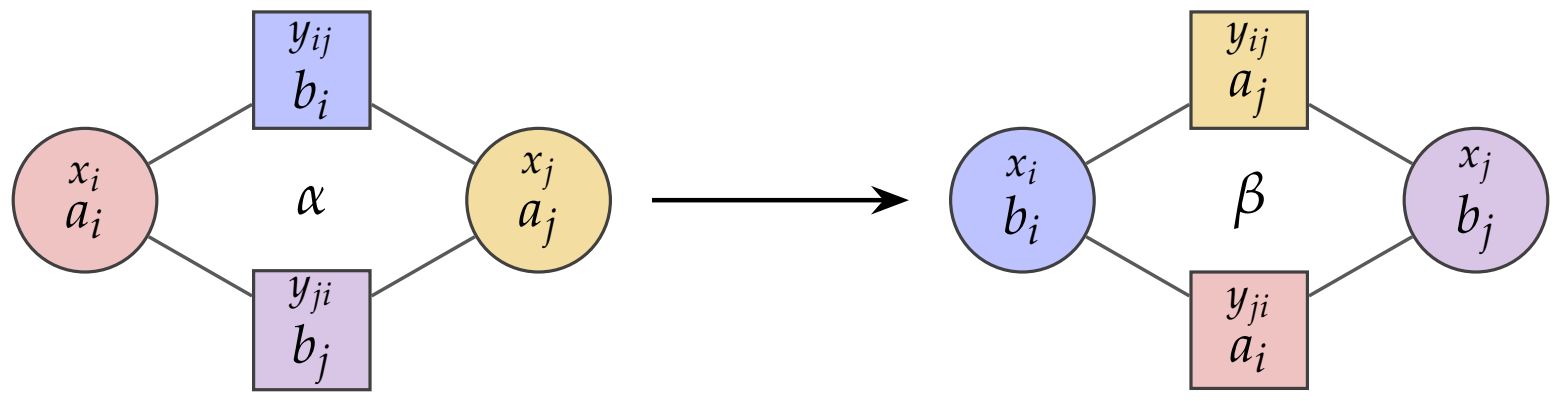
\begin{figure}[h]
\centering
\begin{tikzpicture}

\begin{scope}[xshift=0cm,yshift=0cm]
  \node at (1.5,0) {$\alpha$};

  \colourvertex{fv1Axi}{(0,0)}{x_i}{a_i}{colourAone}
  \colourvertex{fv1Alast}{(3,0)}{x_j}{a_j}{colourAtwo}

  \begin{scope}[square vertices]
    \colourvertex{fv1Ayij}{(1.5,0.86)}{y_{ij}}{b_i}{colourBone}
    \colourvertex{fv1Ayji}{(1.5,-0.86)}{y_{ji}}{b_j}{colourBtwo}
  \end{scope}

  \draw[graph edge] (fv1Axi)--(fv1Ayij);
  \draw[graph edge] (fv1Ayij)--(fv1Alast);
  \draw[graph edge] (fv1Axi)--(fv1Ayji);
  \draw[graph edge] (fv1Ayji)--(fv1Alast);
\end{scope}

\begin{scope}[xshift=6.2cm,yshift=0cm]
  \node at (1.5,0) {$\beta$};

  \colourvertex{fv1Bxi}{(0,0)}{x_i}{b_i}{colourBone}
  \colourvertex{fv1Blast}{(3,0)}{x_j}{b_j}{colourBtwo}

  \begin{scope}[square vertices]
    \colourvertex{fv1Byij}{(1.5,0.86)}{y_{ij}}{a_j}{colourAtwo}
    \colourvertex{fv1Byji}{(1.5,-0.86)}{y_{ji}}{a_i}{colourAone}
  \end{scope}

  \draw[graph edge] (fv1Bxi)--(fv1Byij);
  \draw[graph edge] (fv1Byij)--(fv1Blast);
  \draw[graph edge] (fv1Bxi)--(fv1Byji);
  \draw[graph edge] (fv1Byji)--(fv1Blast);
\end{scope}

\draw[recolouring arrow] (3.75,0)--(5.45,0);

\end{tikzpicture}
\caption{The two colourings $\alpha$ (left) and $\beta$ (right) restricted to the induced subgraph $G_{\tpar}[\{x_i,y_{ij},y_{ji},x_j\}] \cong G_2$.
}
\label{fig:four-vertices}
\end{figure}

Both colourings $\alpha$ and $\beta$ are proper, and they differ at every
vertex. Therefore every vertex will have to be recoloured at least once, requiring at least $t^2$ recolourings. Consider any shortest recolouring sequence from $\alpha$ to $\beta$, and write its length as
$\tpar^2+r$, where $r$ counts recolourings beyond the first at each vertex.
Call a vertex \emph{repeated} if it is recoloured at least twice.
\newpage

\begin{claim}\label{cl:twice-recoloured}
If $x_i$ is recoloured exactly twice, then it has a repeated
neighbour in $Y$.
\end{claim}

\begin{claimproof}
Suppose $x_i$ is recoloured exactly twice, following
$a_i\to c\to b_i$ for some $c\notin\{a_i,b_i\}$.
If $c=b_j$ for some $j\ne i$, then $y_{ji}$ must leave $b_j$
before $x_i$ leaves $a_i$. Its first recolouring therefore cannot
be directly to its target colour $a_i$, so $y_{ji}$ is repeated.
If $c=a_j$ for some $j\ne i$, suppose that $y_{ij}$ is not repeated.
Its unique recolouring must then be from $b_i$ to $a_j$.
This cannot occur before $x_i$ gets colour $a_j$, as it would block
that recolouring. It also cannot occur while $x_i$ has colour $a_j$,
as this would violate properness. It also cannot occur after $x_i$
is coloured $b_i$, as $y_{ij}$ would still have colour $b_i$ at that time.
Thus $y_{ij}$ is repeated.
\end{claimproof}

\begin{claim}\label{cl:once-recoloured}
If $i\ne j$ and $x_i,x_j$ are both recoloured exactly once,
then at least one of $y_{ij},y_{ji}$ is repeated.
\end{claim}

\begin{claimproof}
Suppose that neither $y_{ij}$ nor $y_{ji}$ is repeated.
Then each is recoloured exactly once.
Writing $u\prec v$ when the unique recolouring of $u$ precedes
the unique recolouring of $v$, properness requires
$x_j\prec y_{ij}\prec x_i$
  and
  $x_i\prec y_{ji}\prec x_j.$
These chronological orders are incompatible.
Hence at least one of $y_{ij},y_{ji}$ is repeated.
\end{claimproof}

Let $X_1,X_2,X_3$ denote the disjoint sets of $X$-vertices that are recoloured once,
twice, and at least three times, respectively, so
$|X_1|+|X_2|+|X_{3}|=\tpar$. Let $y$ be the number of repeated $Y$-vertices. Then
$$r\ge |X_2|+2|X_3|+y.$$
By Claim~\ref{cl:once-recoloured}, there are at least $\binom{|X_1|}{2}$ repeated
$Y$-vertices whose endpoints are both in $X_1$.

Disjoint from that, by Claim~\ref{cl:twice-recoloured}, there are at least $|X_2|/2$ repeated $Y$-vertices adjacent to some vertex of $X_2$, where the division by $2$ is because an $Y$-vertex could be adjacent to two vertices from $X_2$. Hence
$$ y\ge\binom{|X_1|}{2}+\frac{|X_2|}{2}.$$
Therefore
$2r\ge3|X_2|+4|X_3|+|X_1|^2 -|X_1|
 =3{\tpar}-4+(|X_1|-2)^2+|X_3|\ge3{\tpar}-4.$
It follows that the diameter of $\C_{2\tpar}(G_{\tpar})$ is at least
$$\dist_{\C_{2\tpar}(G_{\tpar})}(\alpha,\beta) =t^2+r \ge \tpar^2+\frac{3{\tpar}-4}{2} = 
n(G_{\tpar})+\frac{3{\mu(G_t)}}{2}-2.$$
\end{proof}

\section{Regular counterexamples}\label{sec:regular}
In this section we prove Theorem~\ref{thm:regular} via the following slight strengthening and its Corollary~\ref{cor:8regular}.
\begin{theorem}\label{thm:regular_precise}
For every even $\Delta$, there exists a $\Delta$-regular graph $H$ such that
$$\diam\C_{\Delta+2}(H)
  \ge n(H)+\mu(H)+\left\lfloor\frac{\Delta-4}{4}\right\rfloor.$$
\end{theorem}

\begin{proof}
Fix an even $\Delta\ge 2$, put $t=\Delta/2+1$ and $m=t(t-1)$,
and take the graph $G_t$ with parts $X,Y$ and colourings $\alpha,\beta$ that were constructed in the proof of Theorem~\ref{thm:lower}. Thus $|Y|=m$, every $x_i$
has degree $\Delta$, and every $y_{ij}$ has degree two.

Let $B$ be obtained from the complete bipartite graph $K_{\Delta,\Delta}$ by deleting a matching
of size $(\Delta-2)/2$. Exactly $\Delta-2$ vertices of $B$ have degree
$\Delta-1$, and the remaining vertices have degree $\Delta$. Moreover, note that
$n(B)=2\Delta$ and $B$ has a perfect matching so that $\mu(B)=\Delta$. For each $i\ne j$, take a disjoint
copy $B_{ij}$ of $B$ and join its $\Delta-2$ degree-$(\Delta-1)$ vertices to $y_{ij}$.
The resulting graph $H_t$ is $\Delta$-regular, and
$$ n(H_t)=t^2+2\Delta m,\qquad \mu(H_t)=t+\Delta m.$$
For the matching-number lower bound, combine a maximum matching of
$G_t$ with a perfect matching in each $B_{ij}$. Conversely, a matching
uses at most one edge between $y_{ij}$ and $B_{ij}$ for each $i\neq j$. If it uses one, then at most
$\Delta-1$ edges remain inside $B_{ij}$, and otherwise at most $\Delta$. Thus any given matching has at most $\Delta$  edges with an endpoint in 
$B_{ij}$, for every $i\neq j$, and at most $\mu(G_t)=t$ further edges lie in $G_t$.

Next we lower bound the diameter of the reconfiguration graph of $H_t$. Extend $\alpha,\beta$ to $H_t$ using the same palette $\{a_1,\ldots, a_t,b_1,\ldots,b_t\}$ of
$2t=\Delta+2$ colours that was used for $G_t$. In each $B_{ij}$, choose two distinct colours
$c_{ij},d_{ij}$ outside $\{b_i,a_j\}$, colour its two parts
$c_{ij},d_{ij}$ in $\alpha$, and interchange these colours in $\beta$.
Both extensions are proper, since
$\alpha(y_{ij})=b_i$ and $\beta(y_{ij})=a_j$.

Consider any recolouring sequence between these extensions. Its
restriction to $G_t$ requires at least
$t^2+\lceil3t/2\rceil-2$ recolourings, by the proof of
Theorem~\ref{thm:lower} and integrality. In each $B_{ij}$, every vertex
must be recoloured, and the colours on each edge of a perfect matching
are interchanged. At least one endpoint of each such edge must therefore
be recoloured twice. Hence each copy requires at least $3\Delta$
recolourings. These vertex sets are disjoint, so
$$\diam\C_{\Delta+2}(H_t)
  \ge t^2+\left\lceil\frac{3t}{2}\right\rceil-2+3\Delta m
=n(H_t)+\mu(H_t)+\left\lceil\frac t2\right\rceil-2=n(H_t)+\mu(H_t)+\left\lfloor\frac{\Delta-4}{4}\right\rfloor.$$
\end{proof}


\begin{corollary}\label{cor:8regular}
There exist infinitely many  $8$-regular graphs $H$ with
$\diam\C_{10}(H)
  \ge n(H)+ \tfrac{166}{165}\, \mu(H).
$
\end{corollary}
\begin{proof}
The construction in the proof of Theorem~\ref{thm:regular_precise} yields an $8$-regular graph $H_5$ with 
matching number $5+8 \cdot 20=165$, such that $\diam \C_{10}(H_5) \ge n(H_5) + \mu(H_5) + 1$. Now take $r$ disjoint copies of $H_5$. The resulting graph $H$ satisfies
$\diam \C_{10}(H) \ge r\cdot (n(H_5) + \mu(H_5)+1)= r\cdot (n(H_5) + \tfrac{166}{165}\,\mu(H_5))=n(H)+ \tfrac{166}{165}\,\mu(H).$
\end{proof}

We remark that the construction and computation of Corollary~\ref{cor:8regular} generalizes to all even $\Delta$.
For every even integer $\Delta\ge4$, there exist infinitely many 
$\Delta$-regular graphs $H$ such that \\
$\diam\C_{\Delta+2}(H) \ge n(H)+
  \left(1+\frac{4\floor{(\Delta-4)/4}}{(\Delta+2)(\Delta^2+2)}
  \right)\mu(H).$ Furthermore, 
 we elected to give the simplest construction and proof for Theorem~\ref{thm:regular_precise}, by using the gadget $B$ obtained from $K_{\Delta,\Delta}$ minus a matching. With other choices it is possible to do slightly better. For instance,  Theorem~\ref{thm:regular_precise} also holds with $B=K_{\Delta+2} - (C_{\Delta-2} \cup2 K_2)$, a complete graph minus the disjoint union of a cycle of length $\Delta-2$ and two edges. Since the resulting $\Delta$-regular graph $H$ is smaller, for $\Delta=8$ the analogue of Corollary~\ref{cor:8regular} then gives a better bound $\diam\C_{10}(H)
  \ge n(H)+ \tfrac{106}{105}\, \mu(H)$.

\FloatBarrier
\section{Bipartite graphs with degree at most two on one side}\label{sec:upper}

In this section we prove~Theorem~\ref{thm:upper}. The lower bound directly follows from \cite[Proposition~3]{CCC}: $\diam\C_k(G)\ge n(G)+\mu(G)$ for every bipartite graph $G$ and every $k\ge3$. Therefore we only need to prove the upper bounds, which we will do in a stronger form to facilitate induction. 
Given a graph $G$ and a list-assignment $L:V(G)\rightarrow 2^{\mathbb{N}}$, an \emph{$L$-colouring} (or \emph{list colouring}) is a proper vertex colouring $c:V(G)\rightarrow\mathbb{N}$ such that $c(v)\in L(v)$ for every vertex $v$. The $L$-colouring reconfiguration graph $\C_L(G)$ of a graph $G$ has as vertices the $L$-colourings of $G$, two of which are adjacent iff they differ on precisely one vertex of $G$. Since $k$-colourings (for $k\in \mathbb{N}$) are special cases of list colourings, the following implies the upper bounds in Theorem~\ref{thm:upper}. 
\newpage
\begin{theorem}\label{thm:upper_stronger}
Let $G$ be bipartite with bipartition $A\cup B$, where $\deg(b)\le2$ for every
$b\in B$. 
\begin{enumerate}
\item If $|L(v)|\ge \deg(v)+2$ for every vertex $v$, then $\diam\C_L(G)\le n(G)+\tfrac32\mu(G)$.
\item If $|L(a)|\ge \deg(a)+3$ for every $a\in A$ and $|L(b)|\ge \deg(b)+2$ for every $b\in B$, then
      $\diam\C_L(G)\le n(G)+\mu(G)$.
\end{enumerate}
\end{theorem}

\begin{proof}
Fix one of the two statements, and let
$\lambda\in\{\tfrac32,1\}$ be its coefficient of $\mu(G)$.
Suppose that it is false, and choose a counterexample of
minimum order, with $L$-colourings $\alpha,\beta$ at distance
greater than $n+\lambda\mu$, where $n=n(G)$ and $\mu=\mu(G)$.
The cases $n\le1$ are immediate.

We repeatedly use the following observation. If some recolourings
of $\alpha$ and $\beta$ make them agree on a set $S$, using
$r$ recolouring steps in total, we may fix and delete $S$, removing its
colours from the lists of its remaining neighbours. The
conditions on the list sizes are preserved, since each remaining
vertex loses at most one colour for each deleted neighbour.
Minimality therefore gives
\begin{equation}\label{eq:induction}
  \dist_{\C_L(G)}(\alpha,\beta)
  \le r+n-|S|+\lambda \cdot \mu(G-S).
\end{equation}

Let $D$ be the digraph with an arc $u\to v$ whenever
$uv\in E(G)$ and $\beta(u)=\alpha(v)$. (Note that an arc $u \to v$ means that the current colour $\alpha(v)$ of $v$ blocks $u$ from being directly recoloured to its target colour $\beta(u)$.)
We claim that $D$ is strongly connected. Otherwise, choose
a nonempty proper subset $X\subset V(G)$ with no outgoing arc. First
recolour $G[X]$ from $\alpha$ to $\beta$ without recolouring any vertex of $G-X$ (this is possible because there is no arc from $X$ to $G-X$).
Subsequently recolour $G-X$ from $\alpha$ to $\beta$.

Since $\mu(G[X])+\mu(G-X)\le\mu(G)$, applying~inequality~(\ref{eq:induction}) to $S=X$ and to $S=G-X$
gives a sequence of length at most $n+\lambda\mu$, a
contradiction. Thus $D$ is strongly connected. In particular,
$G$ is connected and $\alpha(v)\ne\beta(v)$ for every vertex.

If $G$ is a tree, then $\diam(\C_L(G))\le n(G)+\mu(G)$ follows from~\cite[Corollary 18]{CCC}.
We may therefore assume that $G$ contains a cycle. In this
case $\mu=|A|$. Since $A$ is one part of the bipartition, the inequality $\mu \le |A|$ is immediate. For the reverse inequality, take a spanning unicyclic subgraph,
orient its cycle cyclically, and orient all remaining edges
away from the cycle. Every vertex has indegree one, while
every vertex of $B$ has outdegree at most one. The edges
directed into $A$ therefore form a matching saturating $A$.

If for some $a\in A$ there is a colour $c\in L(a)$ absent from both
$\alpha(N(a))$ and $\beta(N(a))$, then recolour vertex $a$ to $c$ with respect to both $\alpha$ and $\beta$. This requires two recolouring steps in total. Since
$\mu(G-a)\le |A|-1=\mu-1$, inequality (\ref{eq:induction}) applied to $S=\{a\}$ yields $\dist_{\C_L(G)}(\alpha,\beta) \le 2+n-1+\lambda \cdot(\mu-1)  \le n+\lambda \cdot\mu,$ a contradiction. Hence
$$L(a)\subseteq\alpha(N(a))\cup\beta(N(a)) \qquad\text{for every }a\in A.$$
Let $f(v)=\deg_D^+(v)+\deg_D^-(v)$ denote the number of arcs incident to $v$. The multiset $M(a):=\alpha(N(a))\cup\beta(N(a))$ contains $2\deg(a)$ colours (counted with multiplicity) appearing in the neighbourhood of $a$. Among them, exactly $f(a)$ have
colour $\alpha(a)$ or $\beta(a)$, and both colours occur
by strong connectivity. So the number of \emph{distinct} colours $|M(a)|$ in $M(a)$ is at most $2 \deg(a) - (\deg_D^+(a)-1) -  (\deg_D^{-}(a)-1) = 2\deg(a)-f(a)+2 $. Thus,
$$ \deg(a)+2 \le |L(a)| \le |M(a)|
  \le 2\deg(a)-f(a)+2,
$$
so $f(a)\le \deg(a)$. Since every vertex of $D$ has positive
indegree and outdegree,
$$
2|B|\le\sum_{b\in B}f(b)=\sum_{a\in A}f(a) \le\sum_{a\in A}\deg(a)=\sum_{b\in B}\deg(b)\le2|B|.$$
All inequalities must be equalities. In particular,
$f(a)=\deg(a)$ and $|L(a)|=\deg(a)+2$ for every $a\in A$.
This contradicts the statement for
$\lambda=1$, proving the second list-colouring statement.

It remains to finish the case $\lambda=\tfrac32$.
The equality conditions imply that every $b\in B$ has
degree two and exactly one incoming and one outgoing arc
in $D$. They also imply that for every vertex $a\in A$ we have $L(a)=M(a)$ (viewed as a set), and every colour in $L(a) \setminus \{\alpha(a),\beta(a)\}$ occurs
exactly once in the multiset $M(a)$.

Since $D$ is strongly connected and $|A|\ge2$, it contains
a directed path $u\to b\to v$ with distinct $u,v\in A$.
Write
$$x=\alpha(b)=\beta(u),
  \qquad
  y=\beta(b)=\alpha(v).$$
The unique incoming and outgoing arcs at $b$ imply that
$$y\notin\{\alpha(u),\beta(u)\},
  \qquad
  x\notin\{\alpha(v),\beta(v)\}.$$
Consequently, $y\in L(u)$ and it has a unique occurrence in the multiset
$M(u)$, namely $\beta(b)=y$. In particular, $y \notin \alpha(N(u))$.
Similarly, $x\in L(v)$ and it has a unique occurrence in the multiset
$M(v)$,  namely $\alpha(b)=x$. In particular, $x \notin \beta(N(v))$. 

We now choose some free colour $z\in L(b)\setminus\{x,y\}$ and perform the following six
recolourings (three starting from $\alpha$ and three starting from $\beta$), in the indicated order:

$$\begin{array}{c|ccc}
  \text{starting colouring}
    & \text{first} & \text{second} & \text{third}\\
  \hline
  \alpha & u\gets y & b\gets z & v\gets x\\
  \beta  & v\gets x & b\gets z & u\gets y
\end{array}
$$
The uniqueness properties just established, in particular $y\notin \alpha(N(u))$ and $x \notin \beta(N(v))$, ensure the recolourings
at $u$ and $v$ are allowed. Immediately before the recolouring of $b$,
its two neighbours both have colour $y$ in the first row
and colour $x$ in the second, so that recolouring is possible as well.
The resulting colourings agree on $S=\{u,b,v\}$, and $ \mu(G-S)\le |A|-2=\mu-2$.
Applying \eqref{eq:induction}, we obtain the final contradiction
$\dist_{\C_L(G)}(\alpha,\beta)
  \le 6+n-3+\tfrac32(\mu-2)
  =n+\tfrac{3}{2}\mu.$
\end{proof}

\section{A sharp upper bound for bipartite graphs}\label{sec:bipartite}
In this section we prove Lemma~\ref{lem:bipartite}. We restate it for convenience.
\begin{lemma}If
$G$ is bipartite and $k\ge\Delta(G)+2$, then $\diam\C_k(G)\le \tfrac{3}{2}\,n(G).$
\end{lemma}
\begin{proof}
Fix two proper $k$-colourings $\alpha,\beta$ and a bipartition
$V(G)=A\sqcup B$. 
Choose any colour $c\in[k]$. For $X\in\{A,B\}$ put
$$m_X(c)=|\alpha^{-1}(c)\cap X|+|\beta^{-1}(c)\cap X|.$$
Starting from either colouring $\gamma\in\{\alpha,\beta\}$, first
recolour every vertex $v$ of $B$ currently coloured $c$ to a different
available colour. This is possible because at most $d_G(v)+1$
colours are forbidden by the current neighbours and the old colour $c$.
Since $B$ is independent, these recolouring steps do not interfere with each other.
Next, recolour every vertex of $A$ not already coloured $c$ to $c$.
This is allowed because $B$ no longer uses $c$, and $A$ is independent.
The number of recolourings starting from $\gamma$ is at most
$$|\gamma^{-1}(c)\cap B|+|A|-|\gamma^{-1}(c)\cap A|.
$$
The two resulting colourings agree on $A$, which is monochromatic
in colour $c$. Join them by recolouring each vertex of $B$ at most once,
while keeping $A$ fixed, and then reverse the sequence constructed from
$\beta$. This gives an $\alpha$--$\beta$ recolouring sequence of length at most
$\ell_A=2|A|+|B|+m_B(c)-m_A(c).$
Interchanging $A$ and $B$ gives another such recolouring sequence of length at most
$\ell_B=|A|+2|B|+m_A(c)-m_B(c).$
Both constructions use the same arbitrarily chosen colour $c$.
Since $\ell_A+\ell_B=3(|A|+|B|)=3n$, we obtain
$\dist_{\C_k(G)}(\alpha,\beta)
  \le \min\{\ell_A,\ell_B\}
  \le \left\lfloor\frac{3n}{2}\right\rfloor.$
\end{proof}

\begin{corollary}\label{cor:matching}
If $G$ is bipartite and has a perfect 
matching, then
for every $k\ge\Delta(G)+2$,
$\diam\C_k(G)=n(G)+\mu(G)=\lfloor \frac{3}{2} n(G)\rfloor$.
\end{corollary}
\begin{proof}
Take the two colourings obtained by interchanging colours $1$ and $2$
on the bipartition classes. Every vertex must be recoloured, and for
each edge of a maximum matching at least one endpoint must be recoloured
twice. This gives the lower bound $n(G)+\mu(G)=\lfloor \frac{3}{2} n(G)\rfloor$ and
Lemma~\ref{lem:bipartite} gives the upper bound.
\end{proof}

\section{A general lower bound}
This section contains the proof of Theorem~\ref{thm:generallowerbound} in a slightly stronger form, namely Theorem~\ref{thm:lower_stronger} below. For $k\ge 2\Delta$, it recovers the known lower bound $n+\mu$ on the diameter.
For $k=\Delta+2 \ge 4$, it implies $\diam \C_{\Delta+2}(G)\ge n+\frac{\lfloor\Delta/2\rfloor+1}{\Delta}\,\mu$, which is strictly larger than $n+ \frac{1}{2} \mu$. 

As in earlier efforts to prove $\ge n+\mu$ for several graph classes~\cite{CCC,Zeven23}, the high-level idea is to choose a maximum matching $M$ and a proper $k$-colouring $\alpha$, and then obtain a new colouring $\beta$ from $\alpha$ by interchanging the colours along each edge $uv$ of $M$, i.e. setting $\beta(u)=\alpha(v)$  and $\beta(v)=\alpha(u)$. Moreover, for every vertex $u$ not incident to $M$, a colour $\beta(u)$ distinct from $\alpha(u)$ is chosen. If the resulting  $\beta$ is also a \emph{proper} $k$-colouring then we are done: recolouring $\alpha$ to $\beta$ requires at least $3$ recolouring steps on each edge of $M$, so in total at least $n+|M|$ recolourings are needed to go from $\alpha$ to $\beta$. The challenge with this approach is that there can be many bad choices of $M$ and $\alpha$, for which $\beta$ is not proper. To avoid that difficulty, in the argument below, we instead interchange only a fraction of the vertices, in exchange for a weaker bound.

\begin{theorem}\label{thm:lower_stronger}
Let $G$ be a graph with $n$ vertices, maximum degree
$\Delta\ge 1$, and matching number $\mu$. For every integer $k\ge\Delta+2$,
 $$\diam \C_k(G)\ge
  n+
    \min\left(1,\frac{\lfloor k/2\rfloor}{\Delta}\right)\mu.$$
\end{theorem}
\begin{proof}
Fix a maximum matching $M$ and let $U=V(M)$ be the set of vertices that are an endpoint of some edge of $M$.
Let $\pi:U\to U$ be the function that interchanges the endpoints of every edge of $M$. That is,  $\pi(u)=v$ and $\pi(v)=u$ for every edge $uv$ of $M$. Define the
supergraph $\widehat{G}$ of $G[U]$ with vertex set $U$, and edges
$E(\widehat{G})=E(G[U])\cup
  \{\pi(x)\pi(y):xy\in E(G[U])\}.$
Note that $xy \in E(\widehat{G})$ if and only if $\pi(x) \pi(y) \in E(\widehat{G})$. Moreover, $\widehat{G}$ has maximum degree at most $2\Delta-1$.

Take a palette of $2\Delta$ colours equipped with a fixed-point-free
involution $\sigma$. That is, $\sigma:[2\Delta] \rightarrow [2\Delta]$ is a function such that $\sigma(c) \neq c$ and $\sigma(\sigma(c))=c$ for every $c \in [2\Delta]$. We construct a proper colouring $\gamma$ of $\widehat{G}$ such
that
\begin{equation} \label{eq:maintainedproperty}\gamma(\pi(u))=\sigma(\gamma(u))\qquad \text{ for every } u\in U.
\end{equation}
Process the edges of $M$ in arbitrary order. When an edge $uv$ of $M$ is reached, each
endpoint has at most $2\Delta-2$ previously coloured neighbours.
Let $C_u,C_v$ be the colours on those previously
coloured neighbours. Because the previously coloured vertices are the endpoints of some submatching of $M$, the symmetry of $\widehat{G}$ and the maintained property~(\ref{eq:maintainedproperty}) on $\gamma$ give
$C_v=\sigma(C_u)$. Since $|C_u|\le 2\Delta-2$, we can choose a colour $c\notin C_u$. Then $\sigma(c) \notin \sigma(C_u)$ because $\sigma$ is an involution, and $c\neq \sigma(c)$ because $\sigma$ has no fixed point. Therefore we properly may extend the colouring to $\{u,v\}$ by setting
$\gamma(u)=c$ and $\gamma(v)=\sigma(c)$.
This completes the construction of $\gamma$.

For each $c\in [2\Delta]$, consider the set $M_c:= \{ u\pi(u) \in M \mid \gamma(u) \in \{c, \sigma(c)\}\}$ of edges in $M$ whose endpoints are coloured $c$ and $\sigma(c)$.
Since $M_c= M_{\sigma(c)}$ for all $c$, this partitions $M$ into $\Delta$ disjoint matchings. Put
$r=\min\left(\Delta,\lfloor k/2\rfloor\right)$, and let $M^*\subseteq M$ be the union of the $r$ largest such matchings. Then
$$|M^*|\ge\left\lceil\frac{r}{\Delta}|M|\right\rceil.
$$

Based on $M^*$ and the colouring $\gamma$, we will now define two proper $k$-colourings. Let $U^*\subseteq U$ be the endpoints of $M^*$. The restriction of $\gamma$ to $U^*$ uses at most
$2r\le k$ colours. Rename these colours into $[k]$, and let $\alpha$ denote the resulting
proper $k$-colouring of $\widehat{G}[U^*]$.  As $\widehat{G}$ is a supergraph of $G[U]$, $\alpha$ is a proper $k$-colouring of $G[U^*]$ as well. On $U^*$, define
$\beta(v)=\alpha(\pi(v))$. Then $\beta$ is a proper $k$-colouring of $G[U^*]$, since $xy\in E(G[U^*])$ implies
$\pi(x)\pi(y)\in E(\widehat{G}[U^*])$. By construction, $\alpha$ and $\beta$ differ at every vertex of $U^*$, and
interchange colours across each edge of $M^*$.

Extend both colourings to $G$, one new vertex at a time. At a new vertex
$v$, choose $\alpha(v)$ to avoid its previously assigned neighbours in
$\alpha$. Then choose $\beta(v)$ to avoid its previously assigned neighbours
in $\beta$ and also to avoid $\alpha(v)$. At most $\Delta+1$ colours are
forbidden in the latter step and $k\ge \Delta+2$, so both choices are possible. The resulting proper $k$-colourings differ at every vertex
of $G$.

Every sequence from $\alpha$ to $\beta$ must change every vertex at least
once. On an edge of $M^*$, the first endpoint to change cannot immediately
take its target colour, since that colour is still present at the other
endpoint. Hence at least one endpoint of each edge of $M^*$ needs to be recoloured at least
twice, so in total at least $n+|M^*|$ recolouring steps are needed.
The bound on $|M^*|$ proves the result.
\end{proof}

\section{Extension to hypergraphs}\label{sec:hypergraph}
We end with a rather different topic, namely an extension to hypergraphs. These results were proved in $2024$. The original plan was to write a separate paper about it. To prevent eternal procrastination, we present it here.
Given a hypergraph $\mathcal{H}=(V,E)$, and a vertex $v\in V$, let $\deg^*(v)$ denote the maximum integer $t$ such that there are $t$ hyperedges $e_1,e_2,\ldots, e_t \in E$ such that $v\in e_i$ for all $i$, and $e_i\cap e_j = \{v\}$ for all $1 \leq i <j \leq t$.  In other words, $\deg^*(v)$ is the maximum number of hyperedges containing $v$ that are pairwise vertex-disjoint outside $v$. We can think of $\deg^*(v)$ as the effective degree of $v$. Let $\Delta^*:=\max_{v\in V} \deg^*(v)$ denote the maximum effective degree.
A (proper) $k$-colouring of $\mathcal{H}$ is a mapping $c:V \rightarrow \{1,2,\ldots ,k\}$ such that no hyperedge $e\in E$ is monochromatic. 
The chromatic number $\chi(\mathcal{H})$ of $\mathcal{H}$ is the smallest integer $k$ admitting a proper $k$-colouring. 
Given an integer $k\geq 0$, the $k$-reconfiguration graph $\C_k(\mathcal{H})$ of $\mathcal{H}$ is the graph whose vertices are the proper $k$-colourings of $\mathcal{H}$, any two of which are adjacent if and only if they differ on precisely one vertex. Note that $\C_k(\mathcal{H})$ is non-empty if and only if $k\geq \chi(\mathcal{H})$. Throughout this section we assume that the hypergraph is loopless: every hyperedge of $\mathcal{H}$ contains at least two vertices, so that $\infty \neq \chi(\mathcal{H})\le \Delta^*(\mathcal{H})+1$. 
We show that for such hypergraphs, the diameter of $\C_k(\mathcal{H})$ is also linear in $n$.

\begin{lemma}\label{lem:hypergraph}
Let $\mathcal{H}$ be a loopless hypergraph on $n$ vertices and let $k\geq \Delta^*(\mathcal{H})+2$. Then
$$n \leq \diam \C_k(\mathcal{H}) \leq 2n-1$$
\end{lemma}
\begin{proof}
We first observe that at most $\deg^*(v)$ colours are forbidden
at any vertex $v$. Indeed, suppose there are $t$ distinct colours $d_1,d_2,\ldots, d_t$ such that assigning colour $d_i$ to $v$ renders some hyperedge $F_i$ monochromatic in $d_i$, for every $i\in \{1,2,\ldots t\}$. Then $F_i\cap F_j=\{v\}$ whenever
$i\neq j$, since any vertex in their intersection other than $v$
would have to receive both colours $d_i$ and $d_j$. Thus
$t\leq\deg^*(v)$.

The lower bound is immediate: take any $k$-colouring $\alpha$ and permute the set of colours to obtain another $k$-colouring $\beta$ such that $\alpha(v)\neq \beta(v)$ for every vertex $v$. Since each vertex needs to be recoloured at least once, at least $n$ recolourings are required.

For the upper bound, consider two proper $k$-colourings
$\alpha$ and $\beta$. We describe a procedure that transforms
$\alpha$ into $\beta$, writing $\gamma$ for the current
colouring. Call a colour $d$ \emph{completed} if
$\beta^{-1}(d)\subseteq\gamma^{-1}(d)$. Initially,
$\gamma=\alpha$.
Unless $\gamma=\beta$, there exists an uncompleted colour $c$
such that $|\gamma^{-1}(c)|\leq|\beta^{-1}(c)|$. Indeed,
every completed colour $d$ satisfies
$|\gamma^{-1}(d)|\geq|\beta^{-1}(d)|$. If every uncompleted
colour satisfied the strict reverse of the desired
inequality, summing over all colours would contradict
$\sum_{d\in[k]}|\gamma^{-1}(d)|
=n=\sum_{d\in[k]}|\beta^{-1}(d)|.$

We now perform one phase of the procedure. First, iteratively
recolour every vertex currently coloured $c$ to an available
colour distinct from $c$. By the initial observation and
$k\geq\deg^*(v)+2$, such a colour always exists. 
Once there are no vertices with colour $c$ anymore, proceed with recolouring every vertex in $\beta^{-1}(c)$ with $c$.
The latter recolourings are allowed in arbitrary order, since $\beta^{-1}(c)$ contains no hyperedge. 
By the choice of $c$, this phase requires at most $|\gamma^{-1}(c)|+|\beta^{-1}(c)|\leq 2|\beta^{-1}(c)|$ recolouring steps, and it
makes $c$ completed.

Every previously completed colour $d$ remains completed.
Indeed, the vertices of $\beta^{-1}(d)$ have colour $d$
throughout the phase: they are neither vertices initially
coloured $c$ nor vertices whose target colour is $c$.
Thus each phase increases the number of completed colours,
and no completed colour is ever selected in a later phase.
Repeating the procedure therefore eventually makes every
colour completed, at which point $\gamma=\beta$. Since each
colour is selected at most once and a phase selecting $c$
uses at most $2|\beta^{-1}(c)|$ recolourings, the total number
of recolourings is at most
$
\sum_{c\in[k]}2|\beta^{-1}(c)|=2n.$ 
If some colour $d$ with non-empty $\beta^{-1}(d)$ is never
selected, its contribution is omitted from the sum, saving at
least two steps. Otherwise, in the final phase when the last colour $d$ is selected,
every vertex outside $\beta^{-1}(d)$ already has its target colour. Then we can recolour $\beta^{-1}(d)$ directly to $d$ in at
most $|\beta^{-1}(d)|$ steps, saving at least one step.
\end{proof}

Since hypergraphs can be used to represent a wide variety of combinatorial structures, it is clear that Lemma~\ref{lem:hypergraph} has many corollaries. As an illustration, we mention two of them. 

A vertex-colouring of a digraph is \emph{acyclic} if there is no monochromatic directed cycle. We are interested in the reconfiguration graph whose vertices are the acyclic colourings, two of which adjacent iff they differ on precisely one vertex. 
Let
$\Delta_{\min}(G)=\max_{v\in V(G)}\min(\deg^-(v),\deg^+(v))$. For oriented graphs, Picasarri-Arrieta~\cite[Corollary~11]{PicasarriArrieta2024}
proved that the reconfiguration graph of acyclic $k$-colourings  has diameter at most $2\Delta_{\min} \cdot n$, provided $k\geq\Delta_{\min}+1$ and $\Delta_{\min}\ge 1$. 
For arbitrary digraphs that bound is false. But with one additional allowed colour, the following gives a stronger bound for arbitrary digraphs.

\begin{corollary}
Let $G$ be a digraph with $n$ vertices and let $k$ be at least 
$\Delta_{\min}(G)+2$.
Then the diameter of the reconfiguration graph of acyclic $k$-colourings of $G$ lies between $n$ and $2n-1$.
\end{corollary}
\begin{proof}
Choose $\mathcal{H}$ the hypergraph with vertex set $V(G)$, and with $S\subseteq V(\mathcal{H})$ a hyperedge if and only if $S$ is the vertex set of a directed cycle in $G$. 
Every collection of $t$ directed cycles that contain a vertex $v$ but are disjoint outside $v$ implies the existence of $t$ disjoint pairs $(e^+,e^-)$,
where $e^{+}$ (respectively $e^{-}$)
is an ingoing (respectively outgoing) edge incident to $v$.
Therefore we have $\deg^*(v)\leq \min(\deg^{-}(v), \deg^+(v))$, for every vertex $v$.
\end{proof}

The \emph{vertex arboricity} $a(G)$ of a graph $G$ is the minimum number of colours in a vertex-colouring such that every monochromatic induced subgraph is a forest. Chartrand, Kronk and Wall~\cite{ChartrandKronkWall1968} showed that $a(G)\leq 
\lfloor \frac{\Delta(G)}{2} \rfloor+1$. This also follows
from Lovász's decomposition theorem~\cite{Lovasz1966}. By taking one colour more, we get the reconfiguration analogue.
\begin{corollary}
Let $G$ be an $n$-vertex graph and let $k\ge \lfloor \frac{\Delta(G)}{2} \rfloor +2$. Then the reconfiguration graph of $k$-colourings of $G$ such that every monochromatic induced subgraph is a forest has diameter between $n$ and $2n-1$.
\end{corollary}
\begin{proof}
Choose the hypergraph on $V(G)$ whose hyperedges are the vertex sets of cycles in $G$. Since for every vertex $v$, there are at most $\lfloor \deg(v)/2 \rfloor$ cycles containing $v$ while being pairwise disjoint outside $v$, it follows that $\deg^*(v)\leq \lfloor \deg(v)/2 \rfloor \leq \lfloor \Delta(G)/2 \rfloor$.
\end{proof}

 Let $k(\mathcal{H})$ be the least integer such that $\diam \C_k(\mathcal{H})$ has the same (finite) value for every $k\geq k(\mathcal{H})$.  We call this eventual value the \emph{stable diameter} of $\mathcal{H}$, and $k(\mathcal{H})$ the \emph{stability threshold}. It is not hard to see that the stability threshold is finite. For graphs we know the stable diameter is $n+\mu$, by Proposition~\ref{thm:literature}.\ref{prop1}. This suggests a wider research program: 
\begin{Question}
 Given a hypergraph, what are its stable diameter and stability threshold? 
\end{Question}

\section{Open problems}

While we have disproved all three variants in Conjecture~\ref{conjectures}, it raises more questions than answers. How does the diameter of the reconfiguration graph behave near the connectivity threshold $\Delta+2$? Is there an elegant function that describes it for all graphs, or is the reality much more jagged and littered with exceptions? Which other graph invariants could play a role? Below, we propose some explicit problems for further research, some of which have been posed before.

\begin{itemize}

\item Echoing the final question from~\cite{CCC}: what is the smallest $\epsilon>0$ such that \\$\diam\C_{\Delta(G)+2}(G) \le n(G)+ \epsilon \cdot \mu(G)$  for every graph $G$? By Theorem~\ref{thm:lower} and Proposition~\ref{thm:literature}, we know that 
$3/2 \le \epsilon \le 2$.

\item  
Conversely, what is the largest $\eta>0$ such that $\diam\C_{\Delta(G)+2}(G) \ge n(G)+ \eta \cdot \mu(G)$ for every graph $G$? By Theorem~\ref{thm:generallowerbound} and Proposition~\ref{thm:literature}, we know that $1/2 \le \eta \le 1$. In~\cite{CCC} it was asked whether $\eta=1$.

\item What is the smallest integer $\Delta$ such that $\diam \C_{\Delta+2}(G)$ is strictly larger than $ n(G) + \mu(G)$ for some graph $G$ with maximum degree $\Delta$? The same can be asked with a restriction to $\Delta$-regular graphs. In both cases we currently know via Theorems~\ref{thm:lower} and~\ref{thm:regular} that $\Delta \in \{4,5,6,7,8\}$. 

\item Given a positive integer  $\Delta$, what is the smallest integer $k_{\Delta}$ such that $\diam \C_{k}(G) = n(G)+ \mu(G)$ for all $k\ge k_{\Delta}$ and all graphs $G$ with maximum degree $\Delta$? In other words, what is the \emph{stability threshold} $k_{\Delta}$ such that the diameter settles on the \emph{stable diameter} $n(G)+ \mu(G)$ for all $k \ge k_{\Delta}$? \\
We presently know that $k_{\Delta}=\Delta+2$ for all $\Delta \le 3$ due to De Meyer~\cite{DeMeyer26} but  $k_{\Delta} \ge \Delta+3$ for all even $\Delta \ge 8$ due to Theorem~\ref{thm:lower}. Moreover, $k_\Delta \le 2\Delta+1$ for all $\Delta$ due to~\cite{CCC}, which is sharp for $\Delta=1$. 
For $\Delta\ge 2$, with some care  the upper bound proof can be adapted to yield $k_{\Delta}\le 2\Delta$. In summary, $\Delta+3 \le k_{\Delta} \le 2\Delta$ for infinitely many $\Delta$.

\item 
Classify the graph classes for which 
the diameter of the $k$-colouring reconfiguration graph is exactly equal to $n+\mu$, for all $k \ge \Delta+2$. Is there a characterization in terms of forbidden induced subgraphs?
If equality fails for a certain class, what is then the optimal bounding function in terms of $n$ and $\mu$? 
Other natural graph classes that come to mind are (outer)planar graphs.
\item Analogous questions can be asked for the list colouring and correspondence colouring variants described in \cite{CCC}. In addition, Conjecture $4$ in~\cite{CCC} is still open, concerning list colouring reconfiguration for graphs with bounded maximum average degree.
\end{itemize}

The main initial motivation for Conjecture~\ref{conjectures} was a weaker one: is $\diam \C_{\Delta(G)+2}(G) \le \frac{3}{2} n(G)$, for every graph $G$? 
At the time we were inspired by an argument for complete graphs of Bonamy and Bousquet~\cite[Lemma 5]{BonamyBousquet18} which was later developed further into the Optimal Renaming Lemma~\cite{CambieCamesCranston26}.
De Meyer~\cite{DeMeyer26} proved it is true for all complete multipartite and all subcubic graphs (compare Proposition~\ref{thm:literature}), and it also holds for all bipartite graphs by Lemma~\ref{lem:bipartite}.
Is it true in general?  If so, it is sharp for more than one reason: diameter $\lfloor 3n/2 \rfloor$ is attained by both $K_n$ and $K_{n/2,n/2}$. In terms of only the number of vertices, the current best upper bound (valid for every value of $n$) is $2n-1$, from \cite[Corollary 9]{CCC}.
Since the refinements in Conjecture~\ref{conjectures} have now been disproved, we emphasize and explicitly state 
our (together with Cambie and Cranston) original motivating conjecture.

\begin{conjecture}\label{conj:oldmotivatingconjecture} Let $G$ be a graph. If $k\ge \Delta(G)+2$, then $\diam\C_k(G)\le \tfrac{3}{2}\,n(G).$ \end{conjecture}

\subsection*{Declaration of AI use}
Yes and no.
\subsection*{Declaration of AI use, expanded}
The paper is a balanced recombination of
human and AI-assisted (partial) proofs, both old and new,
that went through dozens of iterations of alternating human and AI-assisted improvements, simplifications and rewriting. The result is a blurry symbiotic mix.

\bibliography{bib}

@article{BartierBousquetFeghaliHeinrichMoorePierron23,
  author  = {Bartier, Valentin and Bousquet, Nicolas and Feghali, Carl and Heinrich, Marc and Moore, Benjamin and Pierron, Th{\'e}o},
  title   = {Recoloring planar graphs of girth at least five},
  journal = {SIAM Journal on Discrete Mathematics},
  volume  = {37},
  number  = {1},
  pages   = {332--350},
  year    = {2023},
  doi     = {10.1137/21M1463598},
  url     = {https://doi.org/10.1137/21M1463598}
}

@article{BartierBousquetHeinrich21,
  author  = {Bartier, Valentin and Bousquet, Nicolas and Heinrich, Marc},
  title   = {Recoloring graphs of treewidth {$2$}},
  journal = {Discrete Mathematics},
  volume  = {344},
  number  = {12},
  pages   = {112553},
  year    = {2021},
  doi     = {10.1016/j.disc.2021.112553},
  url     = {https://doi.org/10.1016/j.disc.2021.112553}
}

@article{BelavadiCameron25,
  author  = {Belavadi, Manoj and Cameron, Kathie},
  title   = {Recoloring some hereditary graph classes},
  journal = {Discrete Applied Mathematics},
  volume  = {361},
  pages   = {389--401},
  year    = {2025},
  doi     = {10.1016/j.dam.2024.10.026},
  url     = {https://doi.org/10.1016/j.dam.2024.10.026}
}

@article{BelavadiCameronHildred25,
  author  = {Belavadi, Manoj and Cameron, Kathie and Hildred, Elias},
  title   = {Frozen colourings in {$2K_2$}-free graphs},
  journal = {Electronic Journal of Combinatorics},
  volume  = {32},
  number  = {2},
  pages   = {P2.29},
  year    = {2025},
  doi     = {10.37236/13669},
  url     = {https://doi.org/10.37236/13669}
}

@article{BelavadiCameronMerkel24,
  author  = {Belavadi, Manoj and Cameron, Kathie and Merkel, Owen},
  title   = {Reconfiguration of vertex colouring and forbidden induced subgraphs},
  journal = {European Journal of Combinatorics},
  volume  = {118},
  pages   = {103908},
  year    = {2024},
  doi     = {10.1016/j.ejc.2023.103908},
  url     = {https://doi.org/10.1016/j.ejc.2023.103908}
}

@article{BelavadiCameronSintiari26,
  author  = {Belavadi, Manoj and Cameron, Kathie and Sintiari, Ni Luh Dewi},
  title   = {Recoloring via modular decomposition},
  journal = {Journal of Graph Theory},
  volume  = {111},
  number  = {4},
  pages   = {113--123},
  year    = {2026},
  doi     = {10.1002/jgt.70005},
  url     = {https://doi.org/10.1002/jgt.70005}
}

@article{BonamyBousquet18,
  author  = {Bonamy, Marthe and Bousquet, Nicolas},
  title   = {Recoloring graphs via tree decompositions},
  journal = {European Journal of Combinatorics},
  volume  = {69},
  pages   = {200--213},
  year    = {2018},
  doi     = {10.1016/j.ejc.2017.10.010},
  url     = {https://doi.org/10.1016/j.ejc.2017.10.010}
}

@article{BonamyJohnsonLignosPatelPaulusma14,
  author  = {Bonamy, Marthe and Johnson, Matthew and Lignos, Ioannis and Patel, Viresh and Paulusma, Dani{\"e}l},
  title   = {Reconfiguration graphs for vertex colourings of chordal and chordal bipartite graphs},
  journal = {Journal of Combinatorial Optimization},
  volume  = {27},
  number  = {1},
  pages   = {132--143},
  year    = {2014},
  doi     = {10.1007/s10878-012-9490-y},
  url     = {https://doi.org/10.1007/s10878-012-9490-y}
}

@article{BousquetFeuilloleyHeinrichRabie25,
  author  = {Bousquet, Nicolas and Feuilloley, Laurent and Heinrich, Marc and Rabie, Mika{\"e}l},
  title   = {Short and local transformations between {$(\Delta+1)$}-colorings},
  journal = {Innovations in Graph Theory},
  volume  = {2},
  pages   = {119--156},
  year    = {2025},
  doi     = {10.5802/igt.8},
  url     = {https://doi.org/10.5802/igt.8}
}

@article{BousquetHeinrich22,
  author  = {Bousquet, Nicolas and Heinrich, Marc},
  title   = {A polynomial version of {Cereceda}'s conjecture},
  journal = {Journal of Combinatorial Theory, Series B},
  volume  = {155},
  pages   = {1--16},
  year    = {2022},
  doi     = {10.1016/j.jctb.2022.01.006},
  url     = {https://doi.org/10.1016/j.jctb.2022.01.006}
}

@article{BousquetPerarnau16,
  author  = {Bousquet, Nicolas and Perarnau, Guillem},
  title   = {Fast recoloring of sparse graphs},
  journal = {European Journal of Combinatorics},
  volume  = {52},
  pages   = {1--11},
  year    = {2016},
  doi     = {10.1016/j.ejc.2015.08.001},
  url     = {https://doi.org/10.1016/j.ejc.2015.08.001}
}

@article{CCC,
  author  = {Cambie, Stijn and {Cames van Batenburg}, Wouter and Cranston, Daniel W.},
  title   = {Optimally reconfiguring list and correspondence colourings},
  journal = {European Journal of Combinatorics},
  volume  = {115},
  pages   = {103798},
  year    = {2024},
  doi     = {10.1016/j.ejc.2023.103798},
  url     = {https://doi.org/10.1016/j.ejc.2023.103798}
}

@article{CambieCamesCranston26,
  author  = {Cambie, Stijn and {Cames van Batenburg}, Wouter and Cranston, Daniel W.},
  title   = {Sharp bounds on lengths of linear recolouring sequences},
  journal = {Electronic Journal of Combinatorics},
  volume  = {33},
  number  = {1},
  pages   = {P1.18},
  year    = {2026},
  doi     = {10.37236/13788},
  url     = {https://doi.org/10.37236/13788}
}

@misc{CambieCamesCranstonHeuvelKang25,
  author  = {Cambie, Stijn and {Cames van Batenburg}, Wouter and Cranston, Daniel W. and van den Heuvel, Jan and Kang, Ross J.},
  title   = {Reconfiguration of list colourings},
  year    = {2025},
  arxiv   = {2505.08020},
  url     = {https://arxiv.org/abs/2505.08020}
}

@phdthesis{cereceda2007diss,
  author  = {Cereceda, Luis},
  title   = {Mixing graph colourings},
  school  = {London School of Economics and Political Science},
  year    = {2007},
  url     = {https://etheses.lse.ac.uk/131/}
}

@article{Cranston22,
  author  = {Cranston, Daniel W.},
  title   = {List-recoloring of sparse graphs},
  journal = {European Journal of Combinatorics},
  volume  = {105},
  pages   = {103562},
  year    = {2022},
  doi     = {10.1016/j.ejc.2022.103562},
  url     = {https://doi.org/10.1016/j.ejc.2022.103562}
}

@article{Cranston25,
  author  = {Cranston, Daniel W.},
  title   = {{$10$}-list recoloring of planar graphs},
  journal = {European Journal of Combinatorics},
  volume  = {130},
  pages   = {104190},
  year    = {2025},
  doi     = {10.1016/j.ejc.2025.104190},
  url     = {https://doi.org/10.1016/j.ejc.2025.104190}
}

@article{DeMeyer26,
  author  = {De Meyer, Lucas},
  title   = {Optimal list recoloring of subcubic graphs and complete multipartite graphs},
  journal = {Discrete Mathematics},
  volume  = {349},
  number  = {8},
  pages   = {115155},
  year    = {2026},
  doi     = {10.1016/j.disc.2026.115155},
  url     = {https://doi.org/10.1016/j.disc.2026.115155}
}

@article{DvorakFeghali21,
  author  = {Dvo{\v{r}}{\'a}k, Zden{\v{e}}k and Feghali, Carl},
  title   = {A {Thomassen}-type method for planar graph recoloring},
  journal = {European Journal of Combinatorics},
  volume  = {95},
  pages   = {103319},
  year    = {2021},
  doi     = {10.1016/j.ejc.2021.103319},
  url     = {https://doi.org/10.1016/j.ejc.2021.103319}
}

@article{Feghali21,
  author  = {Feghali, Carl},
  title   = {Reconfiguring colorings of graphs with bounded maximum average degree},
  journal = {Journal of Combinatorial Theory, Series B},
  volume  = {147},
  pages   = {133--138},
  year    = {2021},
  doi     = {10.1016/j.jctb.2020.11.001},
  url     = {https://doi.org/10.1016/j.jctb.2020.11.001}
}

@article{FeghaliJohnsonPaulusma16,
  author  = {Feghali, Carl and Johnson, Matthew and Paulusma, Dani{\"e}l},
  title   = {A reconfigurations analogue of {Brooks'} theorem and its consequences},
  journal = {Journal of Graph Theory},
  volume  = {83},
  number  = {4},
  pages   = {340--358},
  year    = {2016},
  doi     = {10.1002/jgt.22000},
  url     = {https://doi.org/10.1002/jgt.22000}
}

@misc{GuLeiLiMaWang26,
  author  = {Gu, Ruijuan and Lei, Hui and Li, Zhaoxiang and Ma, Yulai and Wang, Susu},
  title   = {Reconfiguration graphs of {$K_{2,3}$}-minor-free graphs},
  year    = {2026},
  arxiv   = {2605.31225},
  url     = {https://arxiv.org/abs/2605.31225}
}

@incollection{Heuvel13,
  author    = {van den Heuvel, Jan},
  title     = {The complexity of change},
  editor    = {Blackburn, Simon R. and Gerke, Stefanie and Wildon, Mark},
  booktitle = {Surveys in Combinatorics 2013},
  series    = {London Mathematical Society Lecture Note Series},
  volume    = {409},
  pages     = {127--160},
  publisher = {Cambridge University Press},
  year      = {2013},
  doi       = {10.1017/CBO9781139506748.005},
  url       = {https://doi.org/10.1017/CBO9781139506748.005}
}

@article{jerrum1995very,
  author  = {Jerrum, Mark},
  title   = {A very simple algorithm for estimating the number of {$k$}-colorings of a low-degree graph},
  journal = {Random Structures \& Algorithms},
  volume  = {7},
  number  = {2},
  pages   = {157--165},
  year    = {1995},
  doi     = {10.1002/rsa.3240070205},
  url     = {https://doi.org/10.1002/rsa.3240070205}
}

@article{LeiMaMiaoShiWang24,
  author  = {Lei, Hui and Ma, Yulai and Miao, Zhengke and Shi, Yongtang and Wang, Susu},
  title   = {Reconfiguration graphs for vertex colorings of {$P_5$}-free graphs},
  journal = {Annals of Applied Mathematics},
  volume  = {40},
  number  = {4},
  pages   = {394--411},
  year    = {2024},
  doi     = {10.4208/aam.OA-2024-0024},
  url     = {https://doi.org/10.4208/aam.OA-2024-0024}
}

@article{Nishimura18,
  author  = {Nishimura, Naomi},
  title   = {Introduction to reconfiguration},
  journal = {Algorithms},
  volume  = {11},
  number  = {4},
  pages   = {52},
  year    = {2018},
  doi     = {10.3390/a11040052},
  url     = {https://doi.org/10.3390/a11040052}
}

@article{WangLu26,
  author  = {Wang, Yichen and Lu, Mei},
  title   = {Linear recoloring diameter of degenerate chordal graphs and bounded treewidth graphs},
  journal = {Discrete Mathematics},
  volume  = {349},
  number  = {7},
  pages   = {115055},
  year    = {2026},
  doi     = {10.1016/j.disc.2026.115055},
  url     = {https://doi.org/10.1016/j.disc.2026.115055}
}

@mastersthesis{Zeven23,
  author  = {Zeven, K. B.},
  title   = {Reconfiguration of graph colorings},
  school  = {Delft University of Technology},
  type    = {Bachelor's thesis},
  year    = {2023},
  url     = {https://resolver.tudelft.nl/uuid:43108b69-8da3-4a86-b3aa-c4091c13139e}
}

@article{DyerFlaxmanFriezeVigoda06,
  author  = {Dyer, Martin and Flaxman, Abraham D. and Frieze, Alan M. and Vigoda, Eric},
  title   = {Randomly coloring sparse random graphs with fewer colors than the maximum degree},
  journal = {Random Structures \& Algorithms},
  volume  = {29},
  number  = {4},
  pages   = {450--465},
  year    = {2006},
  doi     = {10.1002/rsa.20129},
  url     = {https://doi.org/10.1002/rsa.20129}
}

@article{BonamyBousquet13,
  author  = {Bonamy, Marthe and Bousquet, Nicolas},
  title   = {Recoloring bounded treewidth graphs},
  journal = {Electronic Notes in Discrete Mathematics},
  volume  = {44},
  pages   = {257--262},
  year    = {2013},
  doi     = {10.1016/j.endm.2013.10.040},
  url     = {https://doi.org/10.1016/j.endm.2013.10.040}
}

@article{BonamyBousquetPerarnau21,
  author  = {Bonamy, Marthe and Bousquet, Nicolas and Perarnau, Guillem},
  title   = {Frozen {$(\Delta+1)$}-colourings of bounded degree graphs},
  journal = {Combinatorics, Probability and Computing},
  volume  = {30},
  number  = {3},
  pages   = {330--343},
  year    = {2021},
  doi     = {10.1017/S0963548320000139},
  url     = {https://doi.org/10.1017/S0963548320000139}
}

@article{CranstonMahmoud24,
  author  = {Cranston, Daniel W. and Mahmoud, Reem},
  title   = {{$5$}-coloring reconfiguration of planar graphs with no short odd cycles},
  journal = {Journal of Graph Theory},
  volume  = {105},
  number  = {4},
  pages   = {670--679},
  year    = {2024},
  doi     = {10.1002/jgt.23064},
  url     = {https://doi.org/10.1002/jgt.23064}
}

@article{DaiLiu26,
  author  = {Dai, Tingting and Liu, Runrun},
  title   = {List recoloring of planar graphs without {$4$}-cycles},
  journal = {Discrete Mathematics},
  volume  = {349},
  number  = {11},
  pages   = {115264},
  year    = {2026},
  doi     = {10.1016/j.disc.2026.115264},
  url     = {https://doi.org/10.1016/j.disc.2026.115264}
}

@misc{ZhengLu26,
  author = {Zheng, Ruilin and Lu, Junying},
  title  = {The diameter of recoloring graphs under a maximum average degree bound},
  year   = {2026},
  arxiv  = {2609.32688},
  url    = {https://arxiv.org/abs/2609.32688}
}

@misc{AsgarliKrehbielMacLeanZaimi25,
  author = {Asgarli, Shamil and Krehbiel, Sara and MacLean, Simon and Zaimi, Gjergji},
  title  = {Extremal diameters of {$3$}-coloring graphs of trees},
  year   = {2025},
  arxiv  = {2512.03789},
  url    = {https://arxiv.org/abs/2512.03789}
}

@article{Lovasz1966,
  author  = {L{\'a}szl{\'o} Lov{\'a}sz},
  title   = {On decomposition of graphs},
  journal = {Studia Scientiarum Mathematicarum Hungarica},
  volume  = {1},
  pages   = {237--238},
  year    = {1966}
}

@article{ChartrandKronkWall1968,
  author  = {Gary Chartrand and Hudson V. Kronk and Curtiss E. Wall},
  title   = {The point-arboricity of a graph},
  journal = {Israel Journal of Mathematics},
  volume  = {6},
  pages   = {169--175},
  year    = {1968},
  doi     = {10.1007/BF02760181}
}

@article{PicasarriArrieta2024,
  author  = {Lucas Picasarri-Arrieta},
  title   = {Strengthening the directed {Brooks}' theorem for oriented graphs and consequences on digraph redicolouring},
  journal = {Journal of Graph Theory},
  volume  = {106},
  number  = {1},
  pages   = {5--22},
  year    = {2024},
  doi     = {10.1002/jgt.23066}
}

\end{document}